\documentclass{amsart}
\usepackage{times}      
\usepackage{latexsym}   
\usepackage{amssymb}    
\usepackage{amsmath}    
\usepackage{amsthm}     
\usepackage{mathrsfs}   
\usepackage{euscript}   
\usepackage{epsfig,citesort}             
\usepackage[all]{xy}    

\begin{document}

\newtheorem{theorem}{Theorem}
\newtheorem{corollary}[theorem]{Corollary}
\newtheorem{definition}{Definition}
\newtheorem{lemma}[theorem]{Lemma}
\newcommand{\dif}{\mathrm{d}}
\newcommand{\me}{\mathrm{e}}
\newcommand{\Ker}{\text{Ker}\,}
\newcommand{\Image}{\text{Im}\,}
\newcommand{\M}{\widetilde{M}}
\newcommand{\minus}{\setminus}
\newcommand{\cross}{\times}
\newcommand{\bndry}{\partial{\M}}
\newcommand{\st}{\mid}
\newcommand{\comp}{\circ}
\newcommand{\norm}{\Arrowvert}
\newcommand{\ra}{\rightarrow}
\newcommand{\R}{\mathbb{R}}
\newcommand{\Q}{\mathbb{Q}}
\newcommand{\C}{\mathbb{C}}
\newcommand{\Hy}{\mathbb{H}}
\newcommand{\Ot}{\widetilde{\Omega}}
\newcommand{\Et}{\widetilde{E}}
\newcommand{\At}{\widetilde{A}}
\newcommand{\Bt}{\widetilde{B}}
\newcommand{\Ob}{\overline{\Omega}}
\renewcommand{\epsilon}{\varepsilon}
\newcommand{\nin}{\not\in}
\renewcommand{\qed}{\square}
\renewcommand{\"}[1]{\mathaccent 127 #1}
\renewcommand{\phi}{\varphi}
\newcommand{\ri}{\underline{\rho}}
\newcommand{\rs}{\overline{\rho}}
\newcommand{\te}{t \epsilon}
\newcommand{\cc}{d_{cc}}
\newcommand{\fcc}{f_{cc}}
\newcommand{\db}{\overline{d}}
\newcommand{\pt}{\|_t}
\newcommand{\Obar}{\overline{\Omega}}
\newcommand{\Hams}{Hamenst$\ddot{\text{a}}$dt's }
\newcommand{\hr}{$h$-rank}
\newcommand{\Ham}{Hamenst$\ddot{\text{a}}$dt }
\renewcommand{\"}[1]{\mathaccent 127 #1}
\renewcommand{\phi}{\varphi}
\newcommand{\bmv}{\partial{\M} \setminus \{v_{-\infty}\}}
\newcommand{\bmvw}{\partial{\M} \setminus \{v_{-\infty},w_{-\infty}\}}
\newcommand{\ghlim}{\text{GH-lim}}
\newcommand{\wlim}{\text{$\omega$-lim}}
\newcommand{\Ha}{\mathscr{H}}
\newcommand{\ip}{<\cdot,\cdot>}
\newcommand{\pf}{\noindent {\em Proof: }}
\newcommand{\df}{\underset{def}{\equiv}} 
\newtheorem{Pro}{Proposition}[section]
\newtheorem{Sub}[Pro]{Sublemma}
\newtheorem{Thm}[Pro]{Theorem}
\newtheorem{MThm}{Theorem}
\renewcommand{\theMThm}{\Alph{MThm}}
\newtheorem{Rem}{Remark}
\newtheorem{Def}[Pro]{Definition}
\newtheorem{Not}{Notation}
\newtheorem{Cor}[Pro]{Corollary}
\newtheorem{Ithm}{Theorem}
\renewcommand{\theIthm}{\arabic{chapter}.\arabic{Ithm}}
\newtheorem{IndThm}{Theorem}
\newtheorem{Lem}[IndThm]{Lemma}
\newtheorem{Idef}{Definition}
\renewcommand{\theIdef}{\arabic{chapter}.\arabic{Idef}}
\newcommand{\hook}{\lfloor}

\title{$C^{1}$ hypersurfaces of the Heisenberg group are N-rectifiable}
\author{Daniel R. Cole}
\address{Rice University, Houston, TX 70005}
\email{dcole@math.rice.edu}
\thanks{The first author is supported by NSF grant DMS-0240058}

\author{Scott D. Pauls}
\address{Dartmouth College, Hanover, NH 03755}
\email{scott.pauls@dartmouth.edu}
\thanks{The second author is partially supported by NSF grant DMS-0306752}
\date{\today}

\begin{abstract}
We show that $C^{1}$ hypersurfaces in the Heisenberg group are
countably $N$-rectifiable.  As a corollary, this shows that all
$C^1_H$ graphs over the xy-plane are countable $N$-rectifiable,
showing the equivalence of this notion of rectifiability with that of
Franchi, Serra Cassano and Serapioni \cite{FSSC2} for such surfaces.
\end{abstract}

\maketitle

\section{Introduction}
In this note, we consider two notions of rectifiability for surfaces
in the Heisenberg group.  First, we review and fix notation.  We denote by $\Hy$ the three dimensional Heisenberg group and let
$\mathfrak{h}$ be its Lie algebra.  Recall that 
\[\mathfrak{h} = span \{X,Y,Z\}\]
where the only nontrivial bracket relation is $[X,Y]=Z$.  In this
paper, we will use an identification with $\R^3$ for computational
purposes.  Let $\{x,y,z\}$ be the standard coordinates on $\R^3$.
Then, 
\[X= \partial_x-\frac{y}{2}\partial_z\]
\[Y= \partial_y + \frac{x}{2} \partial_z\]
\[ Z = \partial_z\]
As the exponential map is a diffeomorphism, we will identify $\Hy$ with
$\R^3$ as folows: the triple $(a,b,c)$ denotes the point
$e^{a\;X+b\;Y+c\;Z}$.  We fix a background Riemannian metric, $g$, on
$\Hy$, which makes $\{X,Y,Z\}$ an orthonormal basis at each point.  

To describe the Carnot-Carath\'eodory metric on $\Hy$, we review some
definitions. 
\begin{Def}  The horizontal bundle of $\Hy$ is $H\Hy = span\{X,Y\}$.
\end{Def}
\begin{Def} Let $\mathscr{A}$ be the space of absolutely continuous
  paths in $\Hy$ so that, where the derivative exists, it lies in
  $H\Hy$.
\end{Def}
Next, we define the standard {\bf Carnot-Carath\'eodory metric} on $\Hy$,
\[\cc(m,n)= \inf_{\gamma \in \mathscr{A}}\left \{\int
<\gamma',\gamma'>^\frac{1}{2} | \gamma(0)=m, \gamma(1)=n\right \}\]
For ease of computation, we will use the {\bf Carnot gauge}.  Recall
that if $C$ is a compact set and $m,n \in C$, there exists constants
$B$ (depending on $C$) so that 
\begin{equation}\label{gauge}
B^{-1} \cc(m,n)\le d(m,n) \le B \cc(m,n)
\end{equation}
where $d(m,n)=||m^{-1}n||$ and 
\[||(a,b,c)||=((a^2+b^2)^2+c^2)^\frac{1}{4}\]
See \cite{Bell:1996} for more details on the relation between the
Carnot gauge and the Carnot-Carath\'eodory metric. We note that left
translation in the group and rotations about the $z$-axis are
isometries in both the Carnot-Carath\'eodory metric and the Carnot gauge.

In \cite{Pauls:nrect}, the second author introduced the notion of $N$-rectifiability:

\begin{definition}\label{nrect}  Let $N'$ be a Carnot group and $N$ be a subgroup of $N'$ with Hausdorff dimension $k$.  A subset $E$ of
  another Carnot Group  $(M,d_M)$ is said to be {\bf $N$-rectifiable} if there exists
  $U$, a positive measure subset of $N$, and a Lipschitz map $f:U \ra
  M$ such that 
  $\Ha^k_M(E \minus f(U))=0$.  We say $E$ is {\bf countably
  $N$-rectifiable} if there exist a countable number of pairs of subsets $U_i$ and
  Lipschitz maps $f_i:U_i \ra M$ with $\Ha^k_M(E \minus \cup_i
f_i(U_i))=0$.
\end{definition}
In this definition $\Ha^k_M$ is the $k$-Hausdorff measure on $M$ computed with respect to the Carnot-Carath\'edory metric $d_M$.  In the same paper, the second author showed several local measure
theoretic approximability properties of $N$-rectifiable sets and
demonstrated that level sets of $C^{1}_H$ functions in a Carnot
group share the same approximability properties.  Moreover, he gave a
condition on tangent cones that ensured $N$-rectifiability of subsets of
Carnot groups.  

In \cite{FSSC2}, Franchi, Serra Cassano and Serapioni introduced a
different notion of rectifiability in the Heisenberg group (later, the same authors extended this notion to two step groups \cite{FSSC4,FSSC3}),
$\mathbb{H}$, based on a notion of horizontal regularity:
\begin{Def}\label{C1H} A function $f: \mathbb{H} \ra \R$ is of class $C^1_H$ if
$Xf$ and $Yf$ exist and are continuous.  A surface, $S$, given as the
level set of a function $f$ is called a $C^1_H$ surface if $f\in C^1_H$.
\end{Def} 
We note that while the classes of $C^1$ and $C^1_H$ surfaces are
distinct, if we restrict our attention to graphs over the xy-plane
(i.e. surfaces given by $t=f(x,y)$), then the classes are the same.
Definition \ref{C1H} gives rise to an alternate definition of rectifiability:
\begin{Def}\label{FSSC}
A Cacciopoli set $E$ is $\mathbb{H}$-rectifiable if 
\[\partial E =  N \cup \bigcup_{i=1}^\infty K_i\]
where $N$ has zero $\Ha^3_{\mathbb{H}}$ measure and each $K_i$ is a
compact subset of a noncharacteristic $C^1_H$ hypersurface. 
\end{Def}

In this note, we present a more concrete class of $N$-rectifiable sets
which are not covered in \cite{Pauls:nrect}.  We investigate the
$N$-rectifiability of $C^{1}$
hypersurfaces in the three dimensional Heisenberg group.  

By explicitly constructing Lipschitz mappings, we show the following theorem:

\begin{IndThm}  \label{main}  Let $S$ be a $C^1$ hypersurface in $\Hy$
  and let $N$ be the subgroup given by $\{0,y,z\} \subset \Hy$.  Then
  $S$ is countably N-rectifiable.
\end{IndThm}

The subgroup $N$ in the theorem is
isomorphic to the metric tangent cone at points on $S$ (see
\cite{Bell:1996} or \cite{Mitchell} for a description of the tangent
cones to Carnot-Carath\'eodory metric spaces).  Thus, this
theorem recovers some of the flavor of its Euclidean counterpart:
patches of $C^1$ surfaces look locally like their tangent
approximations.  As pointed out in \cite{Pauls:nrect}, the
Euclidean trick of using the projection is insufficient for the Carnot
setting as the inverse of a projection is not a Lipschitz map.  

As an immediate corollary of to the main theorem in \cite{FSSC2} and the observation after defintion
\ref{FSSC} above, we have:
\begin{Cor}  If $S$ is a $C^1_H$ graph over the xy-plane in
  $\mathbb{H}$ then $S$ is countably $N$-rectifiable.  Moreover, any
  $N$-rectifiable graph over the xy-plane is $\mathbb{H}$-rectifiable.
\end{Cor}  

We briefly remark that this investigation is somewhat different than
that of Serra Cassano and Kirchheim (\cite{SCK}) where those authors
construct Euclidean $\frac{1}{2}$-H$\"\text{o}$lder parameterizations of
$C^1_H$ surfaces and show that the exponent cannot, in general, be
improved.  In our treatment, we insist that the domain be equipped
with a degenerate metric which is inherited from the
Carnot-Carath\'eodory metric on $\mathbb{H}$.

To describe the
mapping geometrically, we make some initial observations.  First, the subgroup $N= \{(0,y,z)\}$ is endowed with a restriction of
the Carnot metric from $\Hy$.  If we consider $N$ abstractly (as
opposed to as a subset of $\Hy$), we will denote the point $(0,y,z)$
by $(y,z)$.  To facilitate checking the Lipschitz condition, we
will use the restriction of the Carnot gauge to compute distances.
Precisely, if $(y_1,z_1),(y_2,z_2) \in N$,
\[ d_N((y_1,z_1),(y_2,z_2)) = \left( (y_1-y_2)^4 + (z_1-z_2)^2\right
)^\frac{1}{4}\]

To understand better the construction of the Lipschitz mapping, we
note that curves of the form $(y,z_0)$, for fixed $z_0$, are
1-Lipschitz images of $\R$ while curves of the form $(y_0,z)$, for
fixed $y_0$, are $\frac{1}{2}$-H$\"\text{o}$lder images of
$\R$. Thus, a necessary condition on the mapping is that these curves
map to curves of the same class. 

Second, if we consider the intersection of the horizontal bundle with
the Riemannian tangent plane at a point on $S$, we have that, generically, there
is a single horizontal line contained in each tangent plane.  To see this, we
compute the Riemannian normal to the surface:
\[ N= (Xf) \; X+ (Yf) \; Y+ (Zf) \; Z\]
We note that the only horizontal vector field (up to a multiple) is
given by
\[ V = -\frac{Yf}{Xf} X + Y \]
We note that 
\[Vf =  -\frac{Yf}{Xf} Xf + Yf = 0\]
and so the vector field $V$ is tangential to the surface $S$ and is
observably horizontal.  The only points where the vector field $V$
would not be well defined (up to a multiple) are the characteristic points,
i.e.\ places where $Xf=Yf=0$.  At points of the charactersitic locus,  we have that the entire tangent space is
horizontal.  
Thus, away from the characteristic locus, there
exists a horizontal field on the surface $S$ and we may consider integral curves of
this vector field.  These curves, as they are horizontal, are
rectifiable curves (in the usual sense) --- locally, they are Lipschitz images of
$\R$.  Other curves on $S$ are $\frac{1}{2}$-H$\"\text{o}$lder images
of $\R$.  From the point of view of understanding the $N$-rectifiability
of the surface, we recall that in
\cite{Balogh}, Z. Balogh showed that the characteristic locus of a $C^1$ hypersurface has
$\Ha^3_{cc}$ measure zero (see also the sharp results of Magnani \cite{Mag:coarea,Mag:coarea2,Mag:blowup}).  As we may ignore sets of measure zero when
discussing $N$-rectifiability, we may ignore the entire characteristic
locus. For the balance of the paper, we consider a noncharacteristic
neighborhood of $S$ and fix a particular choice of integral curves of
$V$.

Putting together these two observations shows how to geometrically
construct the mapping.  On a neighborhood with no characteristic
points, consider a smooth curve $\phi \subset S$,
which is transverse to the horizontal curves everywhere
(i.e.\ $\phi'$ is not a horizontal vector).  Let $\theta_{\phi(s)}(r)$
be the integral curve of the horizontal vector field passing through
$\phi(s)$.  Our candidate mapping is
\begin{equation*}
\begin{split}
\Psi:N &\ra S \subset\Hy \\
 (y,z) &\ra \theta_{\phi(z)}(y)
\end{split}
\end{equation*}

In practice, we will identify a specific curve $\phi$, the
intersection of the surface with the plane $y=0$.  The bulk of the paper is showing that this map is locally Lipschitz.  

\section{Proof of theorem \ref{main}}
To fix notation, we let $f : \mathbb{H} \rightarrow \mathbb{R}$ be a
$C^1$ function and $S$ be the surface given by the level set $f=0$.
Using a suitable left translation, we may assume that $f(\mathbf{0}) = 0$.
Moreover, we assume that the origin is not a characteristic point,
i.e.\ $(Xf,Yf)(\mathbf{0}) \neq \mathbf{0}$.  Again, by composing with suitable
isometries of the Carnot-Carath\'eodory metric, we may assume that
$Xf(\mathbf{0}) > 0$ and $Yf(\mathbf{0}) = 0$.  By the continuity of
$Xf$, there
exists an open neighborhood $U_1$ of $\mathbf{0}$ such that for all $q \in U_1$
we have that $Xf(q) > 0$.

Since $U_1$ is open, there exists an $a > 0$ such that the region \[ \left\{
(x,y,z) \in \mathbb{H} \,\bigg| \, -a \leq x \leq a, \, -a \leq y \leq a, \, -a
- \frac {xy} {2} \leq z \leq a - \frac {xy} {2} \right\} \] is contained in
$U_1$.  Call this region $C_1$.  As $C_1$ is compact, we have that
\begin{itemize}\label{constants}
\item there exists a constant $K > 0$ such that for all $q \in C_1$, $Xf(q)
\geq K$.
\item there exists a constant $L > 0$ such that for all $q \in C_1$,
$|\partial_x f(q)| \leq L$, $|\partial_y f (q)| \leq L$, and $|\partial_z f (q)| \leq L$.
\item there exists a constant $M \geq 0$ such that for all $q \in C_1$, $\left|
\frac{\displaystyle Yf(q)} {\displaystyle Xf(q)} \right| \leq M$.
\end{itemize}
Since on the plane $y = 0$ we have that $X = \partial_x$, if follows
that $K\leq L$.

\noindent
{\bf Claim 1:}  There exists a continuous curve $\phi(z)$
parameterizing the intersection of\\ $f(x,0,z)=0$ and $C_1$ with $z \in
\left [ \frac{-K}{L}a,\frac{K}{L}a \right ]$.

We first note that if $L|z| \leq K|x|$ and $x \geq 0$, then
$f(x, 0, z) \geq 0$:
\begin{align*} f(x, 0, z) &= f(x,0,0) + \int_0^z \partial_z f(x,0,t) \,
  \dif t \geq
f(x,0,0) - L|z| \geq f(x,0,0) - K|x| \\ &= f(0,0,0) + \int_0^x \partial_x f(t,0,0) \,
\dif t - Kx = f(0,0,0) + \int_0^x Xf \, \dif t - Kx \\ &\geq f(0,0,0) + Kx - Kx
= f(0,0,0) = 0; \end{align*} and if $L|z| \leq K|x|$ and $x \leq 0$, then $f(x,
0, z) \leq 0$: \begin{align*} f(x, 0, z) &= f(x,0,0) + \int_0^z \partial_z f(x,0,t)
\, \dif t \leq f(x,0,0) + L|z| \leq f(x,0,0) + K|x| \\ &= f(0,0,0) + \int_0^x
\partial_x f(t,0,0) \, \dif t - Kx = f(0,0,0) + \int_0^x Xf \, \dif t - Kx \\ &\leq
f(0,0,0) + Kx - Kx = f(0,0,0) = 0. \end{align*}

Since $\partial_x f = Xf \geq K$ on the intersection of $C_1$ with the plane $y
= 0$, we have that for all $z \in [-a,a]$ there exists at most one $x \in
[-a,a]$ such that $f(x, 0, z) = 0$.  Fix a value of $z \in \left[ -\frac {K}
{L} a, \frac {K} {L} a \right]$.  Note that then $-a \leq -\frac {L} {K} |z|$
and $\frac {L} {K} |z| \leq a$. The map $x \mapsto f(x, 0, z)$ on the domain
$[-a,a]$ is continuous, and we know from the inequalities above that $f\left(
-\frac {L} {K} |z|, 0, z \right) \leq 0$ and $f\left( \frac {L} {K} |z|, 0, z
\right) \geq 0$.  Thus by the Intermediate Value Theorem there exists $x \in
\left[ -\frac {L} {K} |z|, \frac {L} {K} |z| \right]$ such that $f(x, 0, z) =
0$.  Denote this point $(x, 0, z)$ by $\phi(z)$.  The continuous curve
$\phi(z)$ then parameterizes the intersection of the level set $f(x, y, z) = 0$
with $C_1$ when $y = 0$ and $z \in \left[ -\frac {K} {L} a, \frac {K} {L} a
\right]$, proving the claim.

Next, we recall the vector field defined in the introduction

\[ V = - \frac {Yf} {Xf} X + Y\] on the region
$C_1$.  Note that $V$ is well-defined and continuous everywhere on
$C_1$ as we have assumed that there are no characteristic points
inside $C_1$.

Let $U_2$ be the (non-empty) interior of $C_1$.  Then for each $q \in U_2$
there exists a maximal integral curve $\theta_q(t)$ of $V$, defined on
some open interval $(\alpha(q), \beta(q))$ containing 0, such that $\theta_q
(0) = q$.

Define the set \[ D = \left\{ (0,y,z) \in \mathbb{H} \,\bigg| \, -\frac {K} {L}
a < z < \frac {K} {L} a, \, \alpha(\phi(z)) < y < \beta(\phi(z)) \right\}
\]  On this set the map $(0,y,z) \mapsto \theta_{\phi(z)} (y)$ is well-defined.
The set $D$ is open and contains $\mathbf{0}$, and so there exists $n \leq
\frac {K} {2L}$ such that the set \[ C = 0 \times [-n,n] \times [-n,n] \] is a
subset of $D$.

We define our map $\Psi : C \rightarrow \mathbb{H}$ by the formula \[
\Psi(0,y,z) = \theta_{\phi(z)} (y). \]  Note that by definition $f(\Psi(0,y,z))
= 0$.  We next prove that this is a Lipschitz map from a compact subset of the level
set $x = 0$ to the level set $f(x,y,z) = 0$. \\

Since there exists $L$ such that $d_N((y_1,z_1),(y_2,z_2)) \ge L
((y_1-y_2)^4+(z_1-z_2)^2)^\frac{1}{4}$, it suffices to show that there exists a constant $A \geq 0$ such that for all
$(0,y_1,z_1)$ and $(0,y_2,z_2)$ in $C$, we have that \[ d_\mathbb{H}
(\Psi(0,y_1,z_1), \Psi(0,y_2,z_2)) \leq A \left( (y_1 - y_2)^4 + (z_1 - z_2)^2
\right)^{\frac {1} {4}}. \] 
We first use the triangle inequality to break up
the left hand side: 

\begin{equation*}
\begin{split}
d_\mathbb{H} (\Psi(0,y_1,z_1), \Psi(0,y_2,z_2)) \leq
&d_\mathbb{H} (\Psi(0,y_1,z_1), \Psi(0,y_2,z_1))\\& + d_\mathbb{H}
(\Psi(0,y_2,z_1), \Psi(0,y_2,z_2))
\end{split}
\end{equation*}
We deal with each term on the right
hand side separately, showing that 
\begin{equation}\label{one}
d_\mathbb{H} (\Psi(0,y_1,z_1), \Psi(0,y_2,z_1)) \le
|y_1-y_2|\sqrt{1+M^2}
\end{equation}
and that there exists a constant $A_2$ so that
\begin{equation}\label{two}
d_\mathbb{H} (\Psi(0,y_2,z_1), \Psi(0,y_2,z_2)) \le
A_2 |z_1-z_2|^\frac{1}{2}
\end{equation}

With these estimates in place, we can show that $\Psi$ is Lipschitz:
Let $A_1$ be the maximum of $A_2$ and $\sqrt{ 1 + M^2}$.  Then, putting these
estimates together, we get that 
\begin{align*} d_\mathbb{H} (\Psi(0,y_1,z_1),
\Psi(0,y_2,z_2)) &\leq d_\mathbb{H} (\Psi(0,y_1,z_1), \Psi(0,y_2,z_1))\\&\;\; +
d_\mathbb{H} (\Psi(0,y_2,z_1), \Psi(0,y_2,z_2)) \\ 
&\leq \left| y_1 - y_2
\right| \sqrt{1 + M^2} + A_2 \left| z_1 - z_2 \right|^{\frac {1} {2}} \\ 
&\leq
A_1 \left( \left| y_1 - y_2 \right| + \left| z_1 - z_2 \right|^{\frac {1} {2}}
\right) \\ 
&\leq A_1 \cdot 2^{\frac {3} {4}} \left( (y_1 - y_2)^4 + (z_1 -
z_2)^2 \right)^{\frac {1} {4}} \\ 
&\leq A_1 \cdot 2^{\frac {3} {4}} \cdot B
\cdot d_\mathbb{H} \left( (0, y_1, z_1), (0, y_2, z_2) \right)  \;\;\;
\text{by \eqref{gauge}}\end{align*}
Setting $A = A_1 \cdot 2^{\frac {3} {4}} \cdot B$ proves that $\Psi$
is a Lipschitz map of $C$ onto a neighborhood of the origin in $S$.
Since this construction works for any noncharacteristic point, we can
cover any compact subset of $S$ by a finite number of such
neighborhoods union a portion of the characteristic locus.  By an exhaustion argument, we see that $S$ is the union of
countably many such neighborhoods and the characteristic locus, which
is $\Ha^3_{cc}$-measure zero.  In other words, $S$ is countably N-rectifiable, proving theorem \ref{main}.

We devote the rest of the paper to proving the estimates \ref{one} and
\ref{two}.

\vspace{.1in}

\noindent
{\bf Claim 2:}  \[d_\mathbb{H} (\Psi(0,y_1,z_1), \Psi(0,y_2,z_1)) \le |y_1-y_2|\sqrt{1+M^2}\] 

\vspace{.1in}
To show claim 2, we note that the sub-Riemannian distance between
$\Psi(0,y_1,z_1)$ and $\Psi(0,y_2,z_1)$ is bounded above by the length of any
horizontal curve connecting these two points.  One such curve is
$\theta_{\phi(z_1)} (y)$ for $y_1 \leq y \leq y_2$.  Write this curve in
component form by \[ \phi_{z_1} (y) = (\gamma_1(y), y, \gamma_3(y)). \]  The
length of this curve is
\[ \text{Length} = \left| \int_{y_1}^{y_2} \sqrt{ 1 + \gamma_1'(y)} \, \dif y
\right|, \]
which, since $\theta_{\phi(z_1)} (y)$ is the flow of $V$, is bounded above by
$|y_1 - y_2| \sqrt{1 + M^2}$.  Thus we have that \[ d_\mathbb{H}
(\Psi(0,y_1,z_1), \Psi(0,y_2,z_1)) \leq |y_1 - y_2| \sqrt{1 + M^2} \]

\vspace{.1in}

\noindent
{\bf Claim 3:}{\em  \; \; There exists a constant $A_2$ so that 
\[d_\mathbb{H} (\Psi(0,y_2,z_1), \Psi(0,y_2,z_2)) \le
A_2 |z_1-z_2|^\frac{1}{2}\]}

\vspace{.1in}

First, we write the curve $\theta_{\phi(z_1)} (y)$ for $0 \leq y
\leq y_2$ in the component form $(\gamma_1 (y), y, \gamma_3(y))$ and the curve
$\theta_{\phi(z_2)} (y)$ under the same bounds in the component form $(\eta_1
(y), y, \eta_3(y))$. In terms of these components, we are trying to find a bound on the quantity \[ \left(
(\gamma_1 (y_2) - \eta_1 (y_2))^4 + \left (\gamma_3 (y_2) - \eta_3
(y_2)+\frac{y_2}{2}(\gamma_1(y_2)  - \eta_1 (y_2)) \right)^2
\right)^{\frac {1} {4}} \]
Clearly, it is sufficient to bound $|\gamma_1 (y_2) - \eta_1 (y_2)|$
and $|\gamma_3 (y_2) - \eta_3(y_2)|$.  To this end, we show an intermediate inequality.
\begin{Lem}
For all $y \in [-n, n]$ 
 \[ | \eta_1 (y) - \gamma_1 (y) | \leq \frac {L} {K}
 \left| \eta_3 (y) - \gamma_3 (y) + \frac {1} {2} y (\eta_1 (y) - \gamma_1(y))
 \right| \]
\end{Lem}

\pf  First, we
multiply the left hand side of this inequality by $\frac {K} {K}$: \[ \left|
\eta_1(y) - \gamma_1(y) \right| = \frac {K} {K} \left| \eta_1(y) - \gamma_1(y)
\right| \]  Next, let $\psi_X (t)$ be the integral curve of the vector field
$X$ starting at $(\gamma_1(y), y, \gamma_3(y))$.  Since $Xf \geq K$ on the set
on which we are working, we have that \[ \left| f(\psi_X (\eta_1(y) -
\gamma_1(y)) - f(\psi_X(0))) \right| \geq K \left| \eta_1(y) - \gamma_1(y)
\right| \]  Applying this inequality and using the fact that $f(\psi_X(0)) =
0$, we have that \[ \frac {K} {K} \left| \eta_1(y) - \gamma_1(y) \right| \leq
\frac {1} {K} \left| f(\psi_X (\eta_1(y) - \gamma_1(y))) \right| \]  The vector
field $X = \left(1, 0, \frac {y} {2} \right)$ is constant along any of its
integral curves, and as such we get that 
\begin{equation*}
\begin{split}
 \psi_X (\eta_1(y) - \gamma_1(y))) &=
\left( \gamma_1(y) + (\eta_1(y) - \gamma_1(y)), y, \gamma_3(y) - \frac {y} {2}
\left( \eta_1(y) - \gamma_1(y) \right) \right)\\&= \left( \eta_1(y), y, 
\gamma_3(y) - \frac {y} {2} \left( \eta_1(y) - \gamma_1(y) \right)
\right)
\end{split}
\end{equation*} 
Thus we now have that \[ \left| \eta_1(y) - \gamma_1(y) \right| \leq \frac {1}
{K} \left| f\left( \eta_1(y), y,  \gamma_3(y) - \frac {y} {2} \left( \eta_1(y) -
\gamma_1(y) \right) \right) \right| \]  To finish this proof, we recall that on
the set on which we are working, $\partial_z f \leq L$. This, along with the
fact that $f\left( \eta_1(y), y, \eta_3(y) \right) = 0$,  implies that 
\begin{equation*}
\begin{split}
\bigg | f \left (\eta_1(y),y,\eta_3(y) \right) &- f \left (\eta_1(y),y,\gamma_3(y)-\frac{y}{2} \left (\eta_1(y) - \gamma_1(t)\right) \right ) \bigg |\\
 &=\bigg| f\left( \eta_1(y), y,  \gamma_3(y) - \frac {y} {2} \left( \eta_1(y) -
 \gamma_1(y) \right) \right) \bigg|=G
\end{split}
\end{equation*}
and so, 
\begin{align*} G &\leq L \left| \eta_3(y) - \left( \gamma_3(y) - \frac
{y} {2} \left( \eta_1(y) - \gamma_1(y) \right) \right) \right| \\ 
 &\leq L \left| \eta_3(y) - \gamma_3(y) +
\frac {y} {2} \left( \eta_1(y) - \gamma_1(y) \right) \right| \\
 &\leq L \left| \eta_3(y) - \gamma_3(y) + \frac {y} {2}
\left( \eta_1(y) - \gamma_1(y) \right) \right| \end{align*}  
We now substitute
the right hand side of this inequality to get our result: \[ \left| \eta_1(y) -
\gamma_1(y) \right| \leq \frac {L} {K} \left| \eta_3(y) - \gamma_3(y) + \frac
{y} {2} (\eta_1(y) - \gamma_1(y)) \right| \]
$\qed$

Second, we rewrite the quantity $\eta_3(y) - \gamma_3(y)$ using the
fact that both $\theta_{\phi(z_1)}(y)$ and $\theta_{\phi(z_2)}(y)$ are horizontal:
\begin{align*} \eta_3 (y) - \gamma_3 (y) &= \eta_3 (0) + \int_0^y \eta'_3 (t) \,
\dif t - \gamma_3 (0) - \int_0^y \gamma'_3 (t) \, \dif t \\ &= z_2 - z_1 +
\int_0^y \eta'_3 (t) \, \dif t - \int_0^y \gamma'_3 (t) \, \dif t \\ &= z_2 -
z_1 + \int_0^y \left( \frac {1} {2} \eta_1 (t) - \frac {t} {2} \eta'_1 (t)
\right) \, \dif t - \int_0^y \left( \frac {1} {2} \gamma_1 (t) - \frac {t} {2}
\gamma'_1 (t) \right) \, \dif t \\ &= z_2 - z_1 + \int_0^y (\eta_1 (t) -
\gamma_1 (t)) \, \dif t - \frac {y} {2} (\eta_1 (y) - \gamma_1 (y)) \bigg]_0^y
\\ &= z_2 - z_1 + \int_0^y (\eta_1 (t) - \gamma_1 (t)) \, \dif t - \frac {y}
{2} (\eta_1 (y) - \gamma_1 (y)). \end{align*}  Substituting this into the
inequality above, we get that \[ \left|\eta_1 (y) - \gamma_1 (y) \right| \leq
\frac {L} {K} \left| z_2 - z_1 + \int_0^y (\eta_1 (t) - \gamma_1 (t)) \, \dif t
\right|. \] Next we break up the right hand side using the triangle inequality:
\[ \left| \eta_1 (y) - \gamma_1 (y) \right| \leq \frac {L} {K} \left| z_2 - z_1
\right| + \frac {L} {K} \left| \int_0^y (\eta_1 (t) - \gamma_1 (t)) \, \dif t
\right|. \]  

Thus, by Gronwall's lemma, 
\begin{equation*}\label{B1}
|\eta_1(y_0)-\gamma_1(y_0)| \le \frac{L}{K}\text{e}^{\frac{1}{2}}|z_2-z_1|
\end{equation*}
Note, we use that $y \le \frac{K}{2L}$.

Using this, we
can get a bound on $|\eta_3 (y) - \gamma_3 (y)|$ as well: \begin{align*}
|\eta_3 (y) - \gamma_3 (y)| &= \left| z_2 - z_1 + \int_0^y (\eta_1 (t) -
\gamma_1 (t)) \, \dif t - \frac {y} {2} (\eta_1 (y) - \gamma_1 (y)) \right| \\
&\leq \left| z_2 - z_1 \right| + \left| \int_0^y (\eta_1 (t) - \gamma_1 (t)) \,
\dif t \right| + \left| \frac {y} {2} (\eta_1 (y) - \gamma_1 (y))
\right| \\
&\leq \left| z_2 - z_1 \right| +\frac{\text{e}^{\frac{1}{2}}}{2} \left| z_2 - z_1
\right| + \frac {\text{e}^{\frac{1}{2}}} {4}
\left| z_2 - z_1 \right|  \;\; \text{by \eqref{B1}}\\ 
 &\leq \left( 1+\frac{3\text{e}^{\frac{1}{2}}}{4}\right)
\left| z_2 - z_1 \right| 
\end{align*}  
Thus we have that 
\begin{equation*}
\begin{split} \bigg( (\gamma_1 (y_2) -
\eta_1 (y_2))^4 &+ (\gamma_3 (y_2) - \eta_3 (y_2))^2 \bigg)^{\frac {1} {4}}
\leq\\ &\left( \left( \frac {\text{e}^{\frac{1}{2}}L} {K} \left| z_2 - z_1 \right| \right)^4 + \left(\left( 1+\frac{3\text{e}^{\frac{1}{2}}}{4}\right)
 \left| z_2 - z_1 \right| \right)^2 \right)^{\frac {1} {4}} \\
&\leq \left( \frac {\text{e}^2L^4} {K^4} \left| z_2 - z_1 \right|^2 + \left( 1+\frac{3\text{e}^{\frac{1}{2}}}{4}\right)^2
\right)^{\frac {1} {4}} \left| z_1 - z_2 \right|^{\frac {1} {2}}
\end{split}
\end{equation*}

The quantity \[ \left( \frac {\text{e}^2 L^4} {K^4} \left| z_2 - z_1 \right|^2 + \left( 1+\frac{3\text{e}^{\frac{1}{2}}}{4}\right)^2\right)^{\frac {1} {4}} \] is bounded on our domain: we will call this
bound $A_2$.  Hence we have that \[ \left( (\gamma_1 (y_2) - \eta_1 (y_2))^4 +
(\gamma_3 (y_2) - \eta_3 (y_2))^2 \right)^{\frac {1} {4}} \leq A_2 \left| z_1 -
z_2 \right|^{\frac {1} {2}} \]
This completes the proof of the claim.


\end{document}